\documentclass[12pt,a4paper,twoside]{article}
\usepackage{amsfonts}
\usepackage{amsthm,amsmath,amssymb}
\usepackage{color}
\def\trait #1 #2 #3 {\vrule width #1pt height #2pt depth #3pt}
\def\fin{
    \trait .3 5 0
    \trait 5 .3 0
    \kern-5pt
    \trait 5 5 -4.7
    \trait 0.3 5 0
\medskip}
\newtheorem{teor}{Theorem}[section]
\newtheorem{defin}[teor]{Definition}
\newtheorem{lemm}[teor]{Lemma}
\newtheorem{osse}[teor]{Remark}
\newtheorem{prop}[teor]{Proposition}
\newtheorem{defi}[teor]{Definition}
\newtheorem{coro}[teor]{Corollary}
\newtheorem{prob}[teor]{Problem}
\newtheorem{hypo}[teor]{Hypothesis}
\newcommand{\bele}{\begin{lemm}\begin{sl}}
\newcommand{\enle}{\end{sl}\end{lemm}}
\newcommand{\bedef}{\begin{defi}\begin{sl}}
\newcommand{\eddef}{\end{sl}\end{defi}}
\newcommand{\bete}{\begin{teor}\begin{sl}}
\newcommand{\ente}{\end{sl}\end{teor}}
\newcommand{\beos}{\begin{osse}\begin{rm}}
\newcommand{\eddos}{\end{rm}\end{osse}}
\newcommand{\bepr}{\begin{prop}\begin{sl}}
\newcommand{\empr}{\end{sl}\end{prop}}
\newcommand{\bepro}{\begin{prob}\begin{rm}}
\newcommand{\empro}{\end{rm}\end{prob}}
\newcommand{\bede}{\begin{defin}\begin{sl}}
\newcommand{\edde}{\end{sl}\end{defin}}
\newcommand{\beco}{\begin{coro}\begin{sl}}
\newcommand{\enco}{\end{sl}\end{coro}}
\newcommand{\behy}{\begin{hypo}\begin{sl}}
\newcommand{\enhy}{\end{sl}\end{hypo}}

\newcommand{\RR}{\mathbb{R}}

\newcommand{\beeq}[1]{\begin{equation}\label{#1}}
\newcommand{\eddeq}{\end{equation}}
\newcommand{\beeqa}[1]{\begin{eqnarray}\label{#1}}
\newcommand{\eddeqa}{\end{eqnarray}}
\newcommand{\beal}[1]{\begin{align}\label{#1}}
\newcommand{\eddal}{\end{align}}
\newcommand{\bespl}[1]{\begin{split}\label{#1}}
\newcommand{\edspl}{\end{split}}
\newcommand{\bega}[1]{\begin{gather}\label{#1}}
\newcommand{\edga}{\end{gather}}
\newcommand{\beeqax}{\begin{eqnarray*}}
\newcommand{\eddeqax}{\end{eqnarray*}}
\def\qed{\ifmmode   \else \leavevmode\unskip\penalty9999 \hbox{}\nobreak\hfill
  \fi
  \quad\hbox{\hskip.5em\vrule width.4em height.6em depth.05em\hskip.1em}}
\def\endproofsym{\qed}

\def\endnobox{\def\endproofsym{}\end{proof}\def\endproofsym{\qed}}

\newcommand{\no}{\nonumber}
\newcommand{\beeqao}{\begin{eqnarray}\no}
\newcommand{\bealo}{\begin{align}\no}
\newcommand{\besplo}{\begin{split}\no}
\newcommand{\begao}{\begin{gather}\no}
\def\trait #1 #2 #3 {\vrule width #1pt height #2pt depth #3pt}
\def\fin{\hfill
    \trait .3 5 0
    \trait 5 .3 0
    \kern-5pt
    \trait 5 5 -4.7
    \trait 0.3 5 0
\medskip}

\newcommand{\vc}[1]{{\bf #1}}
\newcommand{\dt}{\partial_t}

\newcommand{\dn}{\partial_{\bf n}}

\newcommand{\itt}{\int_0^t}

\newcommand{\io}{\int_\Omega}

\newcommand{\bd}{\boldsymbol{d}}

 \DeclareMathOperator{\dive}{div}

\let\TeXchi\chi
\def\chi{{\setbox0 \hbox{\mathsurround0pt
$\TeXchi$}\hbox{\raise\dp0 \copy0 }}}

\newcommand{\bu}{{\bf u}}
\newcommand{\vr}{\varrho}
\newcommand{\Grad}{\nabla_x}
\newcommand{\tn}[1]{\mbox {\F #1}}

\font\F=msbm10   
scaled 1000

\newcommand{\db}{{\boldsymbol{ d}}}

\newcommand{\ub}{{\bf u}}
\newcommand{\vb}{{\bf v}}

\newcommand{\teta}{\vartheta}

\def\fine{\hfill\kern4pt \vrule height4pt depth0pt width4pt }

\def\dive{\mbox{\rm div\,}}

scaled \magstep1 

   \numberwithin{equation}{section}
\begin{document}

\title{On a non-isothermal model\\ for nematic liquid crystals}

\author{
Eduard Feireisl\thanks{The work of E.F. was supported by
Grant 201/09/0917 of GA \v CR in the framework of research
programmes supported by AV\v CR Institutional Research Plan
AV0Z10190503}\\
Institute of Mathematics of the Czech Academy of Sciences\\
\v Zitn\' a 25, 115 67 Praha 1, Czech Republic\\
feireisl@math.cas.cz
\and
Elisabetta Rocca\thanks{The work of E.R. was
partially supported by the Ne\v cas Center for Mathematical
Modeling sponsored by M\v SMT}\\
Mathematical Department, University of Milan\\
Via Saldini 50, 20133 Milano, Italy
\\
elisabetta.rocca@unimi.it
\and
Giulio Schimperna \\
Mathematical Department, University of Pavia\\
Via Ferrata 1, 27100 Pavia,  Italy
\\
giusch04@unipv.it}

\maketitle

\begin{abstract}
\noindent
A model describing the evolution of a liquid crystal substance in the nematic phase
is investigated in terms of three basic state variables: the {\it absolute temperature} $\teta$,
the {\it velocity field} $\ub$, and the {\it director field} $\bd$, representing
preferred orientation of molecules in a neighborhood of any point
of a reference domain.
The time evolution of the velocity field
is governed by the incompressible Navier-Stokes
system, with a non-isotropic stress tensor depending on the gradients of the
velocity and of the director field $\bd$, where the transport (viscosity) coefficients
vary with temperature. The dynamics of $\bd$ is described by means
of a parabolic equation of Ginzburg-Landau type, with a suitable
penalization term to relax the constraint $|\bd | = 1$. The system
is supplemented by a heat equation, where the heat flux is given
by a variant of Fourier's law, depending also on the
director field $\bd$. The proposed model is shown compatible with
\emph{First and Second laws} of thermodynamics, and the
existence of global-in-time weak solutions for
the resulting PDE system is established, without any essential restriction on the size of the data.
\end{abstract}

\section{Introduction}
\label{sec:intro}

Liquid crystals are materials with rheological properties between
a conventional liquid and a solid crystal, where large elongate
molecules give rise to a preferred orientation. Many different
types of liquid crystals phases have been observed in practical
experiments, distinguished by their characteristic optical
properties (such as birefringence). When viewed under a microscope
with a polarized light source, different liquid crystal phases
will appear to have distinct textures. The contrasting areas in
the textures correspond to domains where the liquid crystals
molecules are oriented in different directions. Within a specific
domain, however, the molecules are well ordered.

Theoretical studies of these types of materials are motivated by real-world applications.
Proper functioning of many practical devices relies on optical properties of certain
liquid crystalline substances in the presence or absence of an
electric field. Typically, a liquid crystal layer sits
between two polarizers that are crossed. The liquid crystal
alignment is chosen so that its relaxed phase is twisted.
This twisted phase reorients light passed through the
first polarizer, allowing its transmission through the second
polarizer (and reflecting back to the observer if a reflector is
provided). The device thus appears transparent. When an electric
field is applied to a liquid crystal layer, the long molecular
axes tend to align parallel to the electric field thus gradually
untwisting in the center of the liquid crystal layer. In such a
state, the liquid crystals molecules do not reorient light, so the
light polarized at the first polarizer is absorbed by the second
polarizer, and the device loses transparency with increasing
voltage. In this way, the electric field can be used to make a
pixel switch between transparent or opaque on command. Color LCD
systems use the same technique, with color filters used to
generate red, green, and blue pixels. Similar principles can be
used to make other liquid crystal based optical devices.

There have been numerous attempts to formulate continuum theories
describing the behavior of liquid crystals flows. We refer to the
seminal papers \cite{eric, le}, where Leslie and Ericksen provide
a mathematical description of various properties of these
materials, and to Lin and Liu \cite{LinLiu} for the first attempt
to analyze the model mathematically. We point out that, to the
present state of knowledge, three main types of liquid crystals
are distinguished, termed \emph{smectic}, \emph{nematic} and
\emph{cholesteric}. The smectic phase forms well-defined layers
that can slide one over another in a manner very similar to that
of a soap. The smectics are ordered along one direction. In the
smectic A phase, the molecules are oriented along the layer
normal, while in the smectic C phase, they are tilted away from
the layer normal. These phases are liquid-like within the layers.
There are many different smectic phases characterized by different
types and degrees of positional and orientational order. The
nematic phase appears to be the most common, where the molecules
do not exhibit any positional order, but they have long-range
orientational order. Thus, the molecules flow and their center of
mass positions are randomly distributed similarly to a liquid, but
they all point in the same direction (within each specific
domain). Most nematics are uniaxial: they possess a preferred axis
that is longer, with the other two being equivalent (so they can be
approximately described as cylinders). Some liquid crystals are
biaxial nematics, meaning that, in addition to orienting their
long axis, they can also orient along a secondary axis. Crystals in
the cholesteric phase exhibit a twisting of the molecules
perpendicular to the director, with the molecular axis parallel to
the director. The main difference between the nematic and
cholesteric phases is that the former is invariant with respect to
certain reflections while the latter is not.

In this paper, we consider the range of temperatures typical for
the nematic phase. As already pointed out, the nematic liquid
crystals are composed of rod-like molecules, with the long axes of
neighboring molecules aligned. Such a kind of anisotropic
structure may be described by means of a dimensionless unit vector
$\bd$, called \emph{director}, that represents the direction of
preferred orientation of molecules in a neighborhood of any point
of a reference domain. In many experiments, the samples of nematic
liquid crystals consist of slowly moving particles, therefore a
relevant approach might be to study the behavior of director field
$\bd$ alone in the absence of velocity. However, the flow velocity
$\bu$ evidently disturbs the alignment of the molecules and also
the converse is true: a change in the alignment will produce a
perturbation of the velocity field $\bu$. Hence, both $\bd$ and
$\bu$ are relevant in the dynamics, and, to a certain extent, also
the changes of the temperature $\teta$ (internal energy). We
introduce a very simple \emph{non-isothermal} model for nematic
liquid crystals in the spirit of the simplified version of the
Leslie-Ericksen model proposed by Lin and Liu \cite{LinLiusimply},
and subsequently studied in \cite{JT, LS}.

In the proposed model, the time evolution of the velocity field
$\bu$ is governed by the standard incompressible Navier-Stokes
system, with a non-isotropic stress tensor depending on $\Grad
\bu$, $\Grad \bd$, where the transport (viscosity) coefficients
vary with temperature. The dynamics of $\bd$ is described by means
of a parabolic equation of Ginzburg-Landau type, with a suitable
penalization term to relax the constraint $|\bd | = 1$. The system
is supplemented by a heat equation, where the heat flux is given
by a variant of Fourier's law in which the dependence on the
director field $\bd$ is taken into account, see Section
\ref{sec:model}. Although such a model may seem rather naive from
the point of view of real-world applications, the present system
of equations captures the essential mathematical features of the
problem, and, last but not least, it is compatible with all
underlying physical principles, in particular with \emph{First and
Second laws} of thermodynamics (as we shall see in Section
\ref{entren}).

Our goal is to show existence of global-in-time weak solutions to
the system, without any essential restriction imposed on the size
of initial data. In order to avoid problems caused by the
interaction of the fluid with a kinematic boundary, we suppose the
latter is impermeable and perfectly smooth imposing the
\emph{complete slip} boundary conditions on the velocity $\bu$.
The existence of weak solutions to the standard incompressible
Navier-Stokes system was established in the celebrated paper by
Leray \cite{leray}. One of the major open problems is to clarify
whether or not the weak solutions also satisfy the corresponding
total energy balance, more precisely, if the kinetic energy of the
system dissipates at the rate given by the viscous stress. Even
the so-called \emph{suitable weak solutions} introduced by
Caffarelli, Kohn and Nirenberg \cite{ckn} allow for an
uncontrolled dissipation of the kinetic energy that may not be
captured by any term appearing in the classical formulation of the
problem. Since the loss of kinetic energy in any energetically
closed system must be compensated by a source term in the internal
energy balance, the above mentioned problem causes unsurmountable
mathematical difficulties whenever the equations for the kinetic
and internal (heat) energy are separated. To avoid this apparent
difficulty, we use the idea proposed in \cite{FM06} replacing the
heat equation by the total energy balance. Of course, the price to
pay is the explicit appearance of the pressure in the total energy
balance that must be handled by refined arguments. Apart from the
fact that the resulting system is mathematically tractable, such
an approach seems much closer to the physical background of the
problem, being an exact formulation of the \emph{First law of
thermodynamics}.

Let us finally mention that, with respect to \cite{FM06}, the main
difficulty here consists in the proof of sufficiently strong
estimates on the director field $\bd$  in order to pass to the
limit in the approximate problem. In particular, the celebrated
Gagliardo-Nirenberg inequality is needed in order to control the
strongly nonlinear terms containing $\Grad\bd$ in both the
momentum equation and the internal energy balance (cf. equations
\eqref{momcons} and \eqref{eqteta} below).

The organization of the paper is as follows. In Section
\ref{sec:model}, we derive the model following the standard
physical principles. The weak formulation, together with the main
result, are stated in Section \ref{sec:mainres}. In Section
\ref{a}, we derive the necessary {\it a priori} estimates and
establish weak compactness of a family of weak solutions subject
to {\it a priori} bounds. Finally, in Section \ref{A}, we
introduce a family of approximate problems, based on Galerkin-type
approximations, and construct a weak solution of the system.

\section{Mathematical model}
\label{sec:model}

Suppose that the fluid occupies a bounded spatial domain $\Omega
\subset \RR^3$, with a sufficiently regular boundary. Let $\vr =
\vr(t,x)$ and $\bu = \bu(t,x)$ denote respectively the mass
density and the velocity in the Eulerian reference system.
Accordingly, the \emph{mass conservation} is expressed by means of
continuity equation
\begin{equation}\label{masscons}
\dt\vr+\dive{(\vr\ub)}=0\,,
\end{equation}
where,
in addition, the standard incompressibility constraint
\begin{equation}\label{incompr}
\dive\ub=0
\end{equation}
is relevant in the context of nematic liquid crystals.

By virtue of \emph{Newton's second law}, the balance of momentum reads
\begin{equation}\label{mombal}
\partial_t (\vr \ub ) + \dive (\vr \ub \otimes \ub )  = \dive \tn{T} + \vr \vc{f} \, ,
\end{equation}
where $\tn{T}$ is the Cauchy stress, and $\vc{f}$ is a given external force.

Motivated by Lin and Liu \cite{LinLiusimply} we consider the stress tensor in the form
\begin{equation}\label{defT}
  \tn{T} = \tn{S} -\vr\lambda(\teta)\left( \Grad \bd \odot \Grad \bd\right) - p\tn{I}\,,
\end{equation}
where $p$ denotes the pressure, and $\tn{S}$ is the conventional Newtonian viscous stress tensor,
\begin{equation}\label{defS}
\tn{S}(\teta, \Grad \ub) = \mu( \teta) \left( \Grad \ub + \Grad^t \ub \right) .
\end{equation}
Note that the transport coefficients $\mu$ and $\lambda$ are
functions of the {\it absolute temperature} $\teta$, see also a strongly related model by
Blesgen \cite{Blesgen}. More specifically, $\mu$ is the viscosity
coefficient assumed always positive, while $\lambda$ denotes the
thermal dilatation coefficient that is an increasing function of
$\teta$.

We assume that the driving force governing the dynamics of the
director $\bd$ is of ``gradient type'' $\partial_{\bd} J$, where
the potential $J$ is given by
\begin{equation}\label{defJ}
J(\teta,\vr,\bd)= W(\bd)+\frac1\teta G(\teta,\vr)\,.
\end{equation}
Here $G$ is a regular function of $\teta$ and $\vr$, and $W$ penalizes the
deviation of the length $|\bd|$ from the value 1. $W$ may be a general function that can be written
as a sum of a convex (possibly non smooth) part, and a smooth, but possibly non-convex one.
A typical example is $W(\bd)=(|\bd|^2-1)^2$. Consequently, $\bd$ satisfies
the following equation
\begin{equation}\label{direq}
\partial_t \bd +  \bu\cdot\Grad \bd+ \partial_{\bd} W(\bd)=\frac{1}{\vr}\dive{(\vr\Grad\bd)}\,.
\end{equation}

Finally, in accordance with the \emph{First law of
thermodynamics}, the internal energy balance reads
\begin{equation}\label{encons}
\dt (\vr e_{\rm int})+\dive(\vr e_{\rm int }\ub)+\dive{\bf q}=\tn{T}:\Grad \ub\,,
\end{equation}
where $e_{\rm int}$ denotes the internal energy density and $\vc{q}$ its flux.
Following Ericksen's model \cite{eric}, the
flux can be taken in the form
\begin{equation} \label{defQ}
{\bf q}=-\kappa (\teta) \Grad\teta-(\kappa_{||}-\kappa_{\bot})(\teta) \db(\db\cdot\Grad\teta),
\end{equation}
where $\kappa,\,\kappa_{||}-\kappa_{\bot}$ are positive functions of the temperature.
Finally, we take $e_{\rm int} = c_v \teta$, where
$c_v > 0$ is the specific heat at constant volume.

Scaling the last equation to have $c_v = 1$, we arrive at the following system:
\begin{align}
\label{eqcont}
& \partial_t \vr+\dive{(\vr\ub)}=0\,,\\
\label{incompressi}
&\dive\ub=0\,,\\
\label{momcons}
& \partial_t(\vr\ub)+\dive(\vr\ub\otimes\ub)+\nabla p=\dive\tn{S}
-\dive{\left(\vr\lambda(\teta)\left(\Grad\bd\odot\Grad\bd\right)\right)} + \vr \vc{f}\,,\\
\label{eqteta}
&\partial_t(\vr\teta)+\dive(\vr\teta \ub)+\dive{\bf q}=\tn{S}:\Grad \ub
-\vr\lambda(\teta)(\Grad\bd\odot\Grad\bd):\Grad \ub\,,\\
\label{eqd}
&\partial_t \bd +\ub\cdot \Grad\bd+ \partial W(\bd)=\frac{1}{\vr}\dive{(\vr\Grad\bd)}\,.
\end{align}

\subsection{Boundary conditions}

Equations (\ref{eqcont} - \ref{eqd}) must be supplemented by a
suitable set of boundary conditions. In order to avoid the occurrence of
boundary layers, we suppose \emph{complete slip} boundary conditions
for the velocity
\begin{equation}\label{slip}
\ub \cdot \vc{n}|_{\partial \Omega} = 0,\ [\tn{T} \vc{n}] \times \vc{n} |_{\partial \Omega} = 0,
\end{equation}
together with the \emph{no-flux} boundary condition for the temperature
\begin{equation} \label{noflux}
\vc{q} \cdot \vc{n}|_{\partial \Omega} = 0,
\end{equation}
and the Neumann boundary condition for the director field
\begin{equation}\label{neumann}
\Grad d_i \cdot \vc{n}|_{\partial \Omega} = 0 \ \mbox{for}\ i=1,2,3.
\end{equation}
The last relation accounts for the fact that there is no contribution to the surface
force $\tn{T} {\bf n}$ from the director $\bd$. This type of boundary conditions not only simplifies the analysis but
it is also suitable for implementation of a numerical scheme (cf.~\cite{LS} for further comments on this topic).

\subsection{Energy, entropy}
\label{entren}

Multiplying momentum equation (\ref{momcons}) by $\vc{u}$ and adding the
resulting expression to (\ref{eqteta}) we deduce the \emph{total energy balance} in the form
\begin{equation}\label{totale}
\partial_t \Big( \vr \left( \frac{1}{2} |\vc{u}|^2 + \teta \right) \Big)
+ \dive \Big( \vr \left( \frac{1}{2} |\vc{u}|^2 + \teta \right) \vc{u} \Big)
+ \dive (p \vc{u}) + \dive \vc{q}
\end{equation}
\[
= \dive (\tn{S} \vc{u}) - \dive{\Big(\vr\lambda(\teta)\left(\Grad\bd\odot\Grad\bd\right) \vc{u} \Big)}
+ \vr \vc{f} \cdot \vc{u}.
\]
Moreover, using the boundary conditions (\ref{slip}), (\ref{noflux})
we may integrate (\ref{totale}) over $\Omega$ to obtain
\begin{equation}\label{enid}
\partial_t\io\left(\frac12\vr|\ub|^2+\vr\teta\right)= \int_{\Omega} \vr \vc{f} \cdot \vc{u},
\end{equation}
in particular, the total energy is a constant of motion as soon as $\vc{f} \equiv 0$.

Let us denote by $\Lambda(\teta)$ a primitive of
$1/\lambda(\teta)$. Testing \eqref{eqteta} by $1/\lambda(\teta)$ and
\eqref{eqd} by $(\dive{(\vr\Grad\bd)}-\vr W'(\bd))$, integrating
the sum of the resulting equations over $\Omega$, and using the boundary
conditions (\ref{neumann}), together with the equation of continuity
\eqref{eqcont}, we get
\begin{align}\no
&\io(\partial_t \bd+\ub\cdot\Grad\bd)\left(\dive{(\vr\Grad\bd)}
-\vr \partial W(\bd)\right)+\dt\io (\vr \Lambda(\teta))
+\io{\bf q}\cdot\Grad\teta\frac{\lambda'(\teta)}{(\lambda(\teta))^2}\\
\label{entropy1}
&=\io\frac{1}{\vr}\left|\dive{(\vr\Grad\bd)}-\vr \partial W(\bd)\right|^2
+\io\frac{1}{\lambda(\teta)}\tn{S}:\Grad \ub
-\io\vr(\Grad\bd\odot\Grad\bd):\Grad \ub\,,
\end{align}
and
\begin{align}\no
&\io(\dt\bd+\ub\cdot\Grad\bd)\left(\dive{(\vr\Grad\bd)}-\vr \partial W(\bd)\right)\\
\no
&= \dt \io \left(-\vr\frac{|\Grad\bd|^2}{2}
-\vr W(\bd)\right)-\io\vr(\Grad\bd\odot\Grad\bd):\Grad \ub\,.
\end{align}
Thus, finally, we arrive at
\begin{equation}\label{entropy2}
\dt \io \vr \left( \Lambda(\teta) -\frac{|\Grad\bd|^2}{2}- \vr W(\bd) \right)
= \io\frac{1}{\vr}\left|\dive{(\vr\Grad\bd)}
-\vr \partial W(\bd)\right|^2
\end{equation}
\[
+\io\frac{1}{\lambda(\teta)}\tn{S}:\Grad \ub-\io{\bf q}\cdot\Grad\teta\frac{\lambda'(\teta)}{(\lambda(\teta))^2},
\]
where
the quantity
$$S=\vr\left(\Lambda(\teta) -|\Grad\bd|^2/2-W(\bd)\right)$$
is the \emph{entropy density} of the system. Accordingly, the
expression under the integral signs on the right-hand side of
(\ref{entropy2}) represents the entropy production. By virtue of
the \emph{Second law of thermodynamics}, the entropy production is
non-negative for any physically admissible process, in particular,
we need  $\lambda' \geq 0$.

\section{Main results}
\label{sec:mainres}

For the sake of simplicity, we restrict ourselves to the case of
constant density, say $\vr \equiv 1$, and $\vc{f} \equiv 0$.
Accordingly, the problem (\ref{incompressi}--\ref{eqd}),
supplemented with the boundary conditions (\ref{slip}),
(\ref{noflux}), and the initial conditions
\begin{equation} \label{initial}
\vc{u}(0, \cdot) = \vc{u}_0, \ \vc{d}(0, \cdot) = \vc{d}_0, \
\teta(0,\cdot) = \teta_0\,,
\end{equation}
reads as follows
\begin{align}
\label{incompressiom}
&\dive\ub=0\,,\\
\label{momconsom}
& \partial_t\ub+\dive(\ub\otimes\ub)+\nabla p=\dive\tn{S}
-\dive{\left(\lambda(\teta)\left(\Grad\bd\odot\Grad\bd\right)\right)}\,,\\
\label{eqtetaom}
&\partial_t\teta+\dive(\teta \ub)+\dive{\bf q}=\tn{S}:\Grad \ub
-\lambda(\teta)(\Grad\bd\odot\Grad\bd):\Grad \ub\,,\\
\label{eqdom}
&\partial_t \bd +\ub\cdot \Grad\bd+ \partial W(\bd)=\Delta\bd\,,
\end{align}
coupled with the boundary conditions (\ref{slip}--\ref{neumann}) and the initial conditions
(\ref{initial}).

To begin, we introduce a weak formulation of (\ref{incompressiom}--\ref{eqdom}) and  formulate
our main result on the existence of global-in-time weak solutions, without any restriction imposed on
the initial data.

\subsection{Weak formulation}
\label{weak}

In the weak formulation, momentum equation (\ref{mombal}), with the incompressibility constraint (\ref{incompr}),
and the boundary conditions (\ref{slip}), are replaced by a family of integral identities
\begin{equation}\label{weak1}
\int_{\Omega} \vc{u}(t, \cdot) \cdot \Grad \varphi = 0 \ \mbox{for a.a.} \ t \in (0,T)
\end{equation}
for any test function $\varphi \in C^{\infty} (\overline{\Omega})$,
\begin{equation} \label{weak2}
\int_0^T \int_{\Omega} \Big( \vc{u} \cdot \partial_t \varphi + \vc{u} \otimes \vc{u} : \Grad \varphi \Big) =
\int_0^T \int_{\Omega} \tn{T} : \Grad \varphi - \int_{\Omega} \vc{u}_0 \cdot \varphi (0, \cdot)
\end{equation}
for any $\varphi \in C^\infty_0 ([0,T) \times \overline{\Omega}; \RR^3)$,
$\varphi \cdot \vc{n}|_{\partial \Omega} = 0$, where the Cauchy stress
tensor $\tn{T}$ is related to the unknowns through the constitutive equation (\ref{defT}).

Equation (\ref{direq}) holds in the strong sense, thanks to the regularity obtained for $\bd$. More specifically, we have
\begin{equation} \label{weak3}
\partial_t \bd + \ub\cdot\Grad \bd + \partial W(\bd) = \Delta \bd\mbox{ a.e. in }(0,T) \times \Omega,
\quad \Grad\bd_i\cdot \vc{n}_{|\partial\Omega}=0\,, \ i=1,2,3\,.
\end{equation}

In the weak formulation the total energy balance (\ref{totale}) is replaced by
\begin{equation} \label{weak4}
\int_0^T \int_{\Omega} \left( \left( \frac{1}{2} |\vc{u}|^2 + \teta \right) \partial_t \varphi +
\left( \frac{1}{2} |\vc{u}|^2 + \teta \right) \vc{u} \cdot \Grad \varphi + \vc{q} \cdot \Grad \varphi \right)
\end{equation}
\[
= \int_0^T \int_{\Omega} \tn{T} \vc{u} \cdot \Grad \varphi - \int_{\Omega}
\left( \frac{1}{2} |\vc{u}_0|^2 + \teta_0 \right)
\varphi (0, \cdot)\,,
\]
for any $\varphi \in C^{\infty}_0 ([0,T) \times \overline{\Omega})$, where $\tn{T}$, $\vc{q}$
obey (\ref{defT}), (\ref{defQ}), and by the {\it entropy inequality}
\begin{equation}\label{weak5}
\partial_t \teta + \dive (\teta\vc{u} ) + \dive
\vc{q} \geq \tn{S}: \Grad \vc{u} - \lambda(\teta)
(\Grad \vc{d} \odot \Grad \vc{d} ) : \Grad \vc{u}\quad\mbox{in }{\cal D}'((0,T)\times\Omega)\,,
\end{equation}
with $\vc{q}$ as in \eqref{defQ} and $\tn{S}$ as in \eqref{defS}.
A {\em weak solution} is a triple $(\ub,\,\bd,\,\teta)$ satisfying
(\ref{weak1}--\ref{weak5}).

\subsection{Main existence theorem}

Before formulating the main result of this paper, let us state the
list of hypotheses imposed on the constitutive functions. We assume
that
\begin{equation}\label{hyp1}
W \in C^2(\RR^3), \quad W \geq 0,\quad \partial W(\vc{d}) \cdot \vc{d} \geq
0 \ \mbox{for all}\ |\vc{d}| \geq D_0
\end{equation}
for a certain $D_0>0$.

In addition, the transport coefficients are continuously
differentiable functions of the absolute temperature satisfying
\begin{equation} \label{hyp2}
0 < \underline{\mu} \leq \mu(\teta) \leq \overline{\mu}, \quad 0 <
\underline{\kappa} \leq \kappa(\teta),\,\,(\kappa_{||} - \kappa_{\bot})
(\teta) \leq \overline{\kappa} \ \mbox{ for all}\ \teta \geq 0
\end{equation}
for suitable constants $\underline{\kappa}$, $\overline{\kappa}$,
$\underline{\mu}$, $\overline{\mu}$.

Finally, let $\lambda\in C^1([0,+\infty))$ be such that
\begin{equation} \label{hyp3}
\lambda'(\teta) \geq 0,\quad \lambda'(0) > 0,\quad \lambda(0) = 0, \quad
\lambda(\teta) \leq \overline{\lambda} \ \mbox{ for all}\ \teta \geq
0
\end{equation}
for a certain $\overline{\lambda} > 0$.

Our main result reads as follows.

\bete \label{theo1} Let $\Omega \subset \RR^3$ be a bounded domain
of class $C^{2 + \nu}$ for some $\nu>0$. Assume that hypotheses (\ref{hyp1} -
\ref{hyp3}) are satisfied. Finally, let the initial data be such
that
\begin{equation}\label{hyp4}
\begin{array}{c}
\vc{u}_0 \in L^2(\Omega; \RR^3),\ \dive \vc{u}_0 = 0,\ \vc{d}_0 \in
L^\infty \cap W^{1,2} (\Omega; \RR^3),\\ \\
 \teta_0 \in
L^1(\Omega),\ {\rm ess} \inf_{\Omega} \teta_0 > 0 .
\end{array}
\end{equation}

Then problem (\ref{weak1}--\ref{weak5}) possesses a weak
solution ($\vc{u}$, $\vc{d}$, $\teta$) in $(0,T) \times \Omega$
belonging to the class
\begin{equation} \label{reg1}
\vc{u} \in L^\infty(0,T; L^2(\Omega; \RR^3)) \cap L^2(0,T;
W^{1,2}(\Omega)),
\end{equation}
\begin{equation} \label{reg2}
\vc{d} \in L^\infty((0,T) \times \Omega; \RR^3) \cap L^\infty(0,T;
W^{1,2}(\Omega; \RR^3))\cap L^2(0,T; W^{2,2}(\Omega;\RR^3)),
\end{equation}
\begin{equation} \label{reg3}
\teta \in L^\infty(0,T; L^1(\Omega)) \cap L^p(0,T;
W^{1,p}(\Omega)),\ 1 \leq p < 5/4, \ \teta > 0 \ \mbox{a.e. in}\
(0,T) \times \Omega,
\end{equation}
with the pressure $p$,
\begin{equation} \label{reg4}
p \in L^{5/3}((0,T) \times \Omega).
\end{equation}
\ente

The rest of the paper is devoted to the proof of Theorem
\ref{theo1}.

\section{A priori bounds}
\label{a}

We establish a number of formal {\it a priori} estimates. These will assume a
rigorous character in the framework of the approximation scheme presented in
Section \ref{A} below.

Combining (\ref{enid}) (multiplied by a positive constant $K>0$)
with (\ref{entropy2}) we obtain the total
dissipation balance in the form
\begin{equation} \label{a1}
\int_{\Omega} \left( \frac{K}{2} |\vc{u}|^2 + ( K\teta -
\Lambda(\teta)) + \frac{|\Grad \vc{d}|^2}{2} + W(\vc{d})
\right)(\tau, \cdot)
\end{equation}
\[
+ \int_0^\tau \int_{\Omega} \left( \left| \Delta \vc{d} - \partial
W(\vc{d}) \right|^2 + \frac{1}{\lambda(\teta)} \tn{S} : \Grad \vc{u}
- \vc{q} \cdot \Grad \teta \frac{ \lambda'(\teta)
}{\lambda^2(\teta)} \right)
\]
\[
\leq \int_{\Omega} \left( \frac{K}{2} |\vc{u}_0|^2 + ( K\teta_0 -
\Lambda(\teta_0)) + \frac{|\Grad \vc{d}_0|^2}{2} + W(\vc{d}_0)
\right).
\]
For $K$ sufficiently large, the terms on the left hand side in \eqref{a1}
turn out to be non-negative, and, in accordance with hypothesis (\ref{hyp4}), the integral on the
right-hand side is bounded; hence we deduce the {\it a priori} bounds
\begin{equation}\label{apr1}
\vc{u} \in L^\infty(0,T; L^2(\Omega; \RR^3)) \cap L^2(0,T;
W^{1,2}(\Omega; \RR^3))\cap L^{10/3}((0,T)\times\Omega;\RR^3),
\end{equation}
\begin{equation}\label{apr2}
\teta, \ \log(\teta) \in L^\infty(0,T; L^1(\Omega)),
\end{equation}
\begin{equation} \label{apr3}
\vc{d} \in L^\infty (0,T; W^{1,2}(\Omega; \RR^3)),
\end{equation}
where we have used (\ref{hyp1} - \ref{hyp3}).

The next step is to take the scalar product of equation
(\ref{weak3}) with $\vc{d}$ yielding
\begin{equation} \label{maxi}
\partial_t |\vc{d}|^2 + \vc{u}\cdot\Grad |\vc{d}|^2 + 2\partial
W(\vc{d}) \cdot \vc{d} = \Delta |\vc{d}|^2 - 2|\Grad \vc{d} |^2.
\end{equation}
By virtue of hypothesis (\ref{hyp1}), we may apply the standard
maximum principle to $|\vc{d}|^2$ to obtain
\begin{equation}\label{apr4}
\vc{d} \in L^\infty ((0,T) \times \Omega; \RR^3).
\end{equation}

Now, going back to (\ref{a1}) and making use of (\ref{apr4}), we get
\begin{equation} \label{apr5}
\vc{d} \in L^2(0,T; W^{2,2}(\Omega; \RR^3)),
\end{equation}
which, together with Gagliardo-Nirenberg interpolation inequality (cf. \cite[p.~125]{nier})
\begin{equation} \label{GL4}
\|\nabla\bd\|_{L^4(\Omega)}\leq
c_1\|\Delta\bd\|_{L^2(\Omega)}^{1/2}\|\bd\|_{L^\infty(\Omega)}^{1/2}+c_2\|\bd\|_{L^\infty(\Omega)}\,,
\end{equation}
gives rise to
\begin{equation}\label{apr6}
\Grad\vc{d} \in L^4((0,T) \times \Omega).
\end{equation}
This estimate turns out to be ``crucial'' in order to obtain a bound for the pressure and, in general,
for the proof of existence of solutions.

Thanks to our choice of the slip boundary conditions (\ref{slip})
for the velocity, the pressure $p$ can be ``computed'' directly from
(\ref{weak3}) as the unique solution of the elliptic problem
\[
\Delta p = \dive \dive \Big( \tn{S} - \lambda(\teta) \Grad \vc{d}
\odot \Grad \vc{d} - \vc{u} \otimes \vc{u} \Big),
\]
supplemented with the boundary condition
\[
\dn p=\left(\dive{\left(\tn{S} - \lambda(\teta) \Grad \vc{d}
\odot \Grad \vc{d} - \vc{u} \otimes \vc{u}\right) }\right)\cdot{\bf n} \ \mbox{ on }\partial\Omega\,.
\]
To be precise, the last two relations have to be interpreted in a ``very weak'' sense. Namely, the pressure
$p$ is determined through a family of integral identities
\begin{equation} \label{press}
\int_\Omega p \Delta \varphi = \int_{\Omega} \Big(\tn{S} -
\lambda(\teta) \Grad \vc{d} \odot \Grad \vc{d} - \vc{u} \otimes
\vc{u} \Big) : \Grad^2 \varphi
\end{equation}
for any test function $\varphi \in C^\infty(\overline{\Omega})$,
$\Grad \varphi \cdot \vc{n}|_{\partial \Omega} = 0$. Consequently,
the bounds established in (\ref{apr1}), (\ref{apr6}) may be used,
together with the standard elliptic regularity results, to conclude
that
\begin{equation}\label{apr7}
p \in L^{5/3}((0,T) \times \Omega).
\end{equation}

Finally, we turn attention to the heat equation (\ref{eqteta}).
Multiplying (\ref{eqteta}) by  $H'(\teta)$ (for a generic $H\in C^2([0,+\infty))$) we deduce its
``renormalized'' form
\begin{equation}\label{renorm}
\partial_t H(\teta) + \dive (H(\teta) \vc{u}) + \dive (H'(\teta)
\vc{q})
\end{equation}
\[
+ H''(\teta) \Big( \kappa(\teta) |\Grad \teta|^2 + (\kappa_{||} -
\kappa_{\bot})(\teta) |\vc{d} \cdot \Grad \teta |^2 \Big)
\]
\[
= H'(\teta) \Big( \tn{S} - \lambda(\teta) \Grad \vc{d} \odot \Grad
\vc{d} \Big) : \Grad \vc{u}\quad\mbox{in }{\cal D}'((0,T)\times\Omega).
\]
The choice $H(\teta) = ( 1 + \teta )^\eta , \ \eta < 1$, in
(\ref{renorm}), together with the uniform bounds obtained in
(\ref{apr1}), (\ref{apr2}), and (\ref{apr6}), yield
\begin{equation} \label{apr8}
\Grad (1 + \teta)^{\nu} \in L^2((0,T) \times \Omega; \RR^3) \
\mbox{for any}\ 0 <\nu < \frac{1}{2}\,.
\end{equation}
Now, we apply an interpolation argument already exploited
 in \cite{BFM}. Using (\ref{apr2}) and \eqref{apr8} and interpolating between
 $\teta\in L^\infty(0,T;L^1(\Omega))$ and $\teta^\nu\in L^1(0,T;L^3(\Omega))$,
 for $\nu\in (0,1]$,
we immediately  get
\begin{equation}\label{apr8bis}
\teta\in {L^q((0,T) \times \Omega)} \ \mbox{for any}\ 1
\leq q < 5/3\,.
\end{equation}
Further, observing that, for all $p\in [1,5/4)$ and $\nu>0$,
\[
\int_{(0,T)\times\Omega}|\nabla\teta|^p\leq
\left(\int_{(0,T)\times\Omega}|\nabla\teta|^2\teta^{\nu-1}\right)^{\frac{p}{2}}
\left(\int_{(0,T)\times\Omega}\teta^{(1-\nu)\frac{p}{2-p}}\right)^{\frac{2-p}{2}}\,,
\]
we conclude from \eqref{apr8} and \eqref{apr8bis} that
\begin{equation} \label{apr9}
\Grad \teta \in L^p((0,T) \times \Omega; \RR^3) \ \mbox{for any}\ 1
\leq p < 5/4.
\end{equation}

The {\it a priori} estimates derived in this section coincide with
the regularity class (\ref{reg1}--\ref{reg4}). Moreover, it can be shown
that the solution set of (\ref{weak1}--\ref{weak4})
is weakly stable (compact) with respect to these
bounds, namely, any sequence of (weak) solutions that complies with
uniform bounds established above has a subsequence that converges to some limit.
Leaving the proof of weak sequential stability to the
interested reader, we pass directly to the proof of Theorem
\ref{theo1} constructing a suitable family of
\emph{approximate} problems whose solutions weakly converges (up to subsequences)
to limit functions which solve the problem in the weak sense specified in Subsection~\ref{weak}.

\section{Approximations}

\label{A}

Solutions to the Navier-Stokes system (\ref{weak1}), (\ref{weak2})
will be constructed by means of the nowadays standard Faedo-Galerkin
approximation scheme, see Temam \cite{temam}. Let
$W^{1,2}_{n,\sigma}(\Omega; \RR^3)$ be the Sobolev space of
solenoidal functions satisfying the impermeability boundary
condition, specifically,
\[
W^{1,2}_{n,\sigma} = \{ \vc{v} \in W^{1,2}(\Omega; \RR^3) \ | \
\dive \vc{v} = 0 \ \mbox{a.e. in}\ \Omega,\ \vc{v} \cdot
\vc{n}|_{\partial \Omega} = 0 \}\,.
\]
Since $\partial \Omega$ is of class $C^{2+ \nu}$, there exists an
orthonormal basis $\{ \vc{v}_n \}_{n=1}^\infty$ of the Hilbert space
$W^{1,2}_{n , \sigma}$ such that $\vc{v}_n \in C^{2 + \nu}$, see
\cite[Theorem 10.13]{FN}. We take $M\leq N$ and denote $X_N = {\rm span}\{ \vc{v}_n
\}_{n=1}^N$. Our strategy is to pass to the limit first for $N\to\infty$ and then
for $M\to\infty$.

The aproximate velocity fields $\vc{u}_{N,M} \in C^1([0,T]; X_N)$ solve
the Faedo-Galerkin system
\begin{equation} \label{approx1}
\frac{{\rm d}}{{\rm d}t} \int_\Omega \vc{u}_{N,M} \cdot \vc{v} =
\int_{\Omega} \vc{u}_{N,M} \otimes [ \vc{u}_{N,M} ]_M : \Grad \vc{v}
\end{equation}
\[
- \int_\Omega \mu(\teta_{N,M}) \Big( \Grad \vc{u}_{N,M} + \Grad^t
\vc{u}_{N,M} \Big): \Grad \vc{v} + \int_\Omega \lambda(\teta_{N,M})
\Grad \vc{d}_{N,M} \odot \Grad \vc{d}_{N,M} : \Grad \vc{v},
\]
\[
\int_\Omega \vc{u}_{N,M}(0, \cdot) \cdot \vc{v} = \int_\Omega
\vc{u}_0 \cdot \vc{v}
\]
for any $\vc{v} \in X_N$. Here, the symbol $[ \vc{v} ]_M$ denotes
the orthogonal projection onto the space ${\rm span}\{V_n
\}_{n=1}^M$.

The functions $\vc{d}_{N,M}$ are determined in terms of
$\vc{u}_{N,M}$ as the unique solution of the parabolic system
\begin{equation} \label{approx2}
\partial_t \vc{d}_{N,M} + [\vc{u}_{N,M}]_M \cdot \Grad \vc{d}_{N,M} +
\partial W(\vc{d}_{N,M}) = \Delta \vc{d}_{N,M},
\end{equation}
supplemented with
\begin{equation} \label{approx3}
\Grad (d_{N,M})_i \cdot\vc{n}|_{\partial \Omega} = 0,\ i = 1,2,3,
\end{equation}
\begin{equation} \label{approx4}
\vc{d}_{N,M}(0, \cdot) = \vc{d}_{0,M},
\end{equation}
where $\vc{d}_{0,M}$ is a suitable smooth approximation of
$\vc{d}_0$.

Next, given $\vc{u}_{N,M}$, $\vc{d}_{N,M}$, the temperature
$\teta_{N,M}$ is evaluated my means of the heat equation (cf. Lady\v
zenskaja et al.~\cite[Chapter V, Theorem 8.1]{lsu})
\begin{align} \label{approx5}
\partial_t \teta_{N,M} + \dive (\teta_{N,M} \vc{u}_{N,M} ) + \dive
\vc{q}_{N,M} = &\ \tn{S}_{N,M}: \Grad \vc{u}_{N,M} \\
\no
&-\left( \lambda(\teta_{N,M}) (\Grad
\vc{d}_{N,M} \odot \Grad \vc{d}_{N,M} )\right) : \Grad \vc{u}_{N,M},
\end{align}
\begin{equation} \label{approx6}
\vc{q}_{N,M} \cdot \vc{n}|_{\partial \Omega} = 0,
\end{equation}
\begin{equation}\label{approx7}
\teta_{N,M} (0, \cdot) = \teta_{0,M},
\end{equation}
where $\tn{S}_{N,M}=\mu( \teta_{N,M}) \left( \Grad \ub_{N,M} + \Grad^t \ub_{N,M} \right)$, and
\[
\vc{q}_{N,M} = - \kappa( \teta_{N,M} ) \Grad \teta_{N,M} -
(\kappa_{||} - \kappa_{\bot} )(\teta_{N,M}) \vc{d}_{N,M} (\vc{d}_{N,M} \cdot
\Grad \teta_{N,M} ).
\]

Finally, the pressure $p_{N,M}$ is determined as the unique solution
of a system of integral identities
\begin{equation} \label{pressap}
\int_\Omega p_{N,M} \Delta \varphi = \int_{\Omega} \Big(\tn{S}_{N,M}
- \lambda(\teta_{N,M}) \Grad \vc{d}_{N,M} \odot \Grad \vc{d}_{N,M} -
\vc{u}_{N,M} \otimes [\vc{u}_{N,M}]_M \Big) : \Grad^2 \varphi
\end{equation}
satisfied
for any test function $\varphi \in C^\infty(\overline{\Omega})$,
$\Grad \varphi \cdot \vc{n}|_{\partial \Omega} = 0$.
In particular, we immediately deduce the estimate
\begin{equation}\no
\|p_{N,M}\|_{L^2((0,T)\times\Omega)}\leq C(M)\,.
\end{equation}

Now, taking $H(\teta)=(1+\teta)^\nu$, with $\nu\in (0,1/2)$, in \eqref{renorm},
we get
$$
\|\dt\teta_{N,M}^\nu\|_{(C^0([0,T];W^{1, r}(\Omega)))^*}\leq C\| \dt\teta_{N,M}^\nu\|_{L^1((0,T)\times\Omega)}\leq C,
$$
where $C$ is a positive constant independent of $N$ and $M$, with $r\in (3,+\infty)$, $\nu\in (0,1/2)$.

Regularizing the convective terms in (\ref{approx1}), (\ref{approx2}) is in the spirit of Leray's original approach
\cite{leray} to the Navier-Stokes system.
As a result, we recover the internal energy {\it equality} at the level of the limit $N \to \infty$. This fact, in turn, enables us to
replace the internal energy equation (\ref{approx5}) by the total energy balance before performing the limit $M \to \infty$.
For fixed $M,N$, problem (\ref{approx1} - \ref{pressap}) can be solved
by means of a simple fixed point argument, exactly as in \cite[Chapter 3]{FN}. Note that all {\it a
priori} bounds derived formally in Section \ref{a} apply to our
approximate problem. Thus given $\vc{u} \in C([0,T]; X_N)$, we find
$\vc{d} = \vc{d}[\vc{u}]$ solving (\ref{approx2} - \ref{approx4}),
and then $\teta = \teta [\vc{u} , \vc{d} ]$ and the pressure $p$
satisfying (\ref{approx5} - \ref{pressap}). Plugging these $\vc{d}$,
$\teta$ in (\ref{approx1}) we may find a new function ${\cal T}[
\vc{u} ]$ defining thus a mapping $\vc{u} \mapsto {\cal T}
[\vc{u}]$. Given the {\it a priori} bounds obtained in Section
\ref{a}, we can easily show that ${\cal T}$ possesses a fixed point
by means of the classical Schauder's argument, at least on a possibly short time interval.
However, using once more the {\it a priori} estimates we easily conclude that the approximate solutions can be extended on any
fixed time interval $[0,T]$ (see \cite[Chapter 6]{FN} for details).


\subsection{Passage to the limit as $N\to\infty$}

Having established the existence of the approximate solutions
$\vc{u}_{N,M}$, $\vc{d}_{N,M}$, $\teta_{N,M}$, and $p_{N,M}$, we
let $N \to \infty$ and use the uniform bounds established in Section
\ref{a} to obtain
\begin{align}
\label{cuN}
&\ub_{N,M}\to\ub_M \ \mbox{ weakly-(*) in } L^\infty(0,T;L^2(\Omega;\RR^3))
\cap L^2(0,T;W^{1,2}(\Omega;\RR^3))\,,\\
\label{cutN}
&\dt\ub_{N,M}\to\dt\ub_M \ \mbox{ weakly in } L^2(0,T;(W^{1,2}(\Omega;\RR^3))^*)\,,\\
\label{cpN}
&p_{N,M}\to p_M \ \mbox{ weakly in } L^{2}((0,T)\times\Omega)\,,\\
\label{ctetaN}
&\teta_{N,M}^\nu\to\teta_M^\nu\ \mbox{ weakly-(*) in } L^2(0,T;W^{1,2}(\Omega))\cap L^\infty(0,T; L^{1/\nu}(\Omega))\,,\\
\label{ctetatN}
&\dt\teta_{N,M}^\nu\to\dt\teta_M^\nu\ \mbox{ weakly-(*) in } (C_0(0,T;W^{1,r}(\Omega)))^*\,,\\
\label{cdN}
&\bd_{N,M}\to \bd_M \ \mbox{ weakly-(*) in } L^\infty(0,T;W^{1,2}(\Omega;\RR^3))
\cap L^4(0,T;W^{2,4}(\Omega;\RR^3))\,,\\
\label{cdtN}
&\dt\bd_{N,M}\to \dt \bd_M \ \mbox{ weakly in } L^4(0,T;L^{4}(\Omega;\RR^3))\,,
\end{align}
for any $\nu\in (0,1/2)$, and $r>3$. Note that at this stage $M$ remains fixed in the convective term
$\ub_{N,M} \otimes[\bu_{N,M}]_M$.

Hence, applying the Aubin-Lions compactness
lemma (cf.~\cite{simon}), we deduce that
\begin{align}\label{custrN}
&\ub_{N,M}\to\ub_M \ \mbox{ strongly in } L^2(0,T;L^2(\Omega;\RR^3))\,,\\
\label{ctetastrN}
&\teta_{N,M}\to \teta \ \mbox{ strongly in } L^p((0,T)\times\Omega)
\end{align}
for any $p\in [1,5/3)$. Moreover, at this level of approximation, the director field $\vc{d}_M$ is regular,
and, in particular, we have
\[
\Grad\bd_{N,M}\to \Grad \bd_M\ \mbox{ strongly in } L^4((0,T)\times\Omega)\,.
\]
Hence, we can perform the limit passage
\[
\lambda (\teta_{N,M}) ( \Grad \vc{d}_{N,M} \odot \Grad \vc{d}_{N,M}
) : \Grad \vc{u}_{N,M}
\]
\[
\to \lambda (\teta_{M}) ( \Grad \vc{d}_{M} \odot \Grad \vc{d}_{M} )
: \Grad \vc{u}_{M} \ \mbox{in, say,}\ L^1((0,T) \times \Omega).
\]

Thus we may infer that the limit quantities $\vc{u}_M$,
$\vc{d}_M$, $\teta_M$, and $p_M$ solve the problem
\begin{equation}\label{aapprox2}
\int_{\Omega} \vc{u}_M(t, \cdot) \cdot \Grad \varphi = 0 \ \mbox{for a.a.} \ t \in (0,T)
\end{equation}
for any test function $\varphi \in C^{\infty} (\overline{\Omega})$,
\begin{equation} \label{aapprox1}
\int_0^T \int_{\Omega} \Big( \vc{u}_M \cdot \partial_t \varphi + \vc{u}_M \otimes [\vc{u}_{M}]_M: \Grad \varphi \Big) =
\int_0^T \int_{\Omega} \tn{T}_M : \Grad \varphi - \int_{\Omega} \vc{u}_0 \cdot \varphi (0, \cdot)
\end{equation}
for any $\varphi \in C^\infty_0 ([0,T) \times \overline{\Omega}; \RR^3)$,
$\varphi \cdot \vc{n}|_{\partial \Omega} = 0$,
\begin{equation} \label{aapprox3}
\partial_t \vc{d}_{M} + [\vc{u}_{M}]_M \cdot \Grad \vc{d}_{M} +
\partial W(\vc{d}_{M}) = \Delta \vc{d}_{M},\quad\mbox{a.e. in }(0,T)\times\Omega,
\end{equation}
supplemented with
\begin{equation} \label{aapprox4}
\Grad (d_{M})_i \cdot\vc{n}|_{\partial \Omega} = 0,\ i = 1,2,3,
\end{equation}
\begin{equation} \label{aapprox5}
\vc{d}_{M}(0, \cdot) = \vc{d}_{0,M},
\end{equation}
and
\begin{equation} \label{aapprox6}
\partial_t \teta_{M} + \dive (\teta_{M} \vc{u}_M ) + \dive
\vc{q}_{M} \geq \tn{S}_{M}: \Grad \vc{u}_M - \lambda(\teta_{M})
(\Grad \vc{d}_{M} \odot \Grad \vc{d}_{M} ) : \Grad \vc{u}_{M}
\end{equation}
in the sense of distributions with non-negative test functions,
\begin{equation} \label{aapprox7}
\vc{q}_{M} \cdot \vc{n}|_{\partial \Omega} = 0,
\end{equation}
\begin{equation}\label{aapprox8}
\teta_{M} (0, \cdot) = \teta_{0,M},
\end{equation}
together with the total energy balance
\begin{equation} \label{aapprox9}
\frac{{\rm d}}{{\rm d}t} \int_{\Omega} \Big( \frac{1}{2}
|\vc{u}_M|^2 + \teta_M \Big) = 0,
\end{equation}
where
\begin{equation}\label{defTM}
  \tn{T}_M = \tn{S}_M -\lambda(\teta_M)\left( \Grad \bd_M \odot \Grad \bd_M\right) - p_M\tn{I}\,,
\end{equation}
and
\begin{equation}\label{defSM}
\tn{S}_M= \mu( \teta_M) \left( \Grad \ub_M + \Grad^t \ub_M \right).
\end{equation}

Moreover, since the convective term  $\vc{u}_M \otimes [\vc{u}_{M}]_M$ is regular (and
consequently $\dt \ub_M\in L^2(0,T;(W^{1,2}(\Omega;\RR^3))^*)$),  we can
take $\bu_M$ as a test function in \eqref{aapprox1} to recover the kinetic energy balance in the form
\begin{align}\no
\|\ub_M(t)\|^2_{L^2(\Omega)}&+\itt\io\mu(\teta_M)|\Grad \vc{u}_{M} + \Grad^t
\vc{u}_{M}|^2=\|\ub_0\|^2_{L^2(\Omega)}\\
\label{compa}
&+2\itt\io\lambda (\teta_{M}) ( \Grad \vc{d}_{M} \odot \Grad \vc{d}_{M} )
: \Grad \vc{u}_{M}\,.
\end{align}
Similarly, taking $\vb =\ub_{N,M}$ in \eqref{approx1} we get
\begin{align}\no
\|\ub_{N,M}(t)\|^2_{L^2(\Omega)}&+\itt\io\mu(\teta_{M,N})|\Grad \vc{u}_{N,M} + \Grad^t
\vc{u}_{N,M}|^2=\|\ub_0\|^2_{L^2(\Omega)}\\
\no
&+2\itt\io\lambda (\teta_{N,M}) ( \Grad \vc{d}_{N,M} \odot \Grad \vc{d}_{N,M} )
: \Grad \vc{u}_{N,M}\,.
\end{align}
Passing to the limit $N\to\infty$ in the last equation, using \eqref{custrN} and \eqref{hyp3},
and comparing the result
with \eqref{compa}, we conclude, by means of  \eqref{hyp2}, that
\begin{equation}\label{custrbisN}
\Grad \ub_{N,M}\to \Grad \ub_M\ \mbox{ strongly in }L^2((0,T)\times\Omega)\,.
\end{equation}
Accordingly, the inequality \eqref{aapprox6} may be replaced by
\begin{equation} \label{aapprox6+}
\partial_t \teta_{M} + \dive (\teta_{M} \vc{u}_M ) + \dive
\vc{q}_{M} = \tn{S}_{M}: \Grad \vc{u}_M - \lambda(\teta_{M})
(\Grad \vc{d}_{M} \odot \Grad \vc{d}_{M} ) : \Grad \vc{u}_{M}\,.
\end{equation}
As a matter of fact, \eqref{aapprox6+} already follows from
(\ref{aapprox6}), (\ref{aapprox9}), and (\ref{compa}).

Moreover, taking in  \eqref{aapprox1} $\ub_M\varphi$ (with $\varphi\in {\cal D}((0,T)\times\Omega)$)
in place of $\varphi$,  we get
\begin{equation}\label{totaleM}
\partial_t  \left( \frac{1}{2} |\ub_M|^2 + \teta_M \right)
+ \dive \Big(\left( \frac{1}{2} |\ub_M|^2 + \teta_M \right) [\ub_M]_M \Big)
+ \dive (p_M \ub_M) + \dive \vc{q}_M
\end{equation}
\[
= \dive (\tn{S}_M \vc{u}_M) - \dive{\Big(\lambda(\teta_M)\left(\Grad\bd_M\odot\Grad\bd_M\right) \vc{u}_M \Big)}
\ \mbox{in }{\cal D}'((0,T)\times\Omega)\,.
\]
This concludes the passage to the limit for $N\to\infty$.


\subsection{Passage to the limit as $M\to\infty$}

Our final goal is to let $M\to\infty$ (\ref{aapprox1}--\ref{aapprox9}).
We notice that the limits in (\ref{cuN}), (\ref{ctetaN}--\ref{ctetatN}), (\ref{custrN}--\ref{ctetastrN})
still hold  when letting  $M\to \infty$.
Moreover, we have
\begin{align}
\label{cdtuM}
&\dt\ub_M\to\dt\ub\ \mbox{ weakly in } L^{5/3}(0,T;W^{-1,5/3}(\Omega;\RR^3))\,,\\
\label{cpM}
&p_M\to p\ \mbox{ weakly in } L^{5/3}((0,T)\times\Omega)\,,\\
\label{cdM}
&\bd_{M}\to \bd\ \mbox{ weakly-(*) in } L^\infty((0,T)\times \Omega; \RR^3)\cap L^\infty(0,T;W^{1,2}(\Omega;\RR^3))\\
\nonumber
&\qquad\qquad\qquad\qquad\qquad\cap L^2(0,T;W^{2,2}(\Omega;\RR^3))\,,\\
\label{cdtM}
&\dt\bd_{M}\to \dt \bd \ \mbox{ weakly in } L^2(0,T;L^{3/2}(\Omega;\RR^3))\,.
\end{align}

Now,
we can easily pass to the limit $M \to \infty$ in (\ref{aapprox2}--\ref{aapprox5}) to recover
(\ref{weak1}--\ref{weak3}).
In addition, by  virtue of \eqref{apr1}, \eqref{apr6}, \eqref{apr8bis}, we get
\begin{align}\nonumber
&\left\{ \left(\frac{|\ub_M|^2}{2}+p_M\right)[\ub_M]_M \right\}_{M > 0} \ \mbox{bounded in } L^{10/9}((0,T)\times\Omega)\,,\\
\nonumber
& \left\{ \teta_M \ub_M \right\}_{M > 0} \ \mbox{bounded in } L^{q}(0,T;L^{q}(\Omega))\ \mbox{for any}\ q\in [1, 10/9)\,,\\
\nonumber
& \left\{ \tn{S}_M \ub_M \right\}_{M > 0} \ \mbox{bounded in } L^{5/4}((0,T)\times\Omega)\,,\\
\nonumber
& \left\{ \lambda(\teta_M) (\Grad\bd_M\odot\Grad\bd_M)\ub_M \right\}_{M > 0}\ \mbox{ bounded in } L^{5/4}(0,T;L^{5/4}(\Omega))\,.
\end{align}
Consequently, we can pass to the limit in \eqref{totaleM} to deduce the desired conclusion \eqref{weak4}.

Finally, as convex functionals are weakly lower semicontinuous, we can see that
\eqref{aapprox6} gives rise to
\begin{equation} \label{entropyineq}
\partial_t \teta + \dive (\teta \vc{u} ) + \dive
\vc{q} \geq \tn{S}: \Grad \vc{u} - \lambda(\teta)
(\Grad \vc{d} \odot \Grad \vc{d} ) : \Grad \vc{u}\ \mbox{ in }{\cal D}'((0,T)\times\Omega)\,.
\end{equation}
This completes the proof of Theorem
\ref{theo1}.

\end{document}